\documentclass{amsart}
\usepackage[T1]{fontenc} 
\usepackage{amssymb}
\usepackage{mathrsfs}
\usepackage{hyperref}
\usepackage{verbatim}
\usepackage{graphicx}
\usepackage{subfigure}
\usepackage{psfrag}
\usepackage{enumitem}

\usepackage[dvipsnames]{xcolor} 

\newtheorem{theorem}{Theorem}[section]

\newtheorem{corollary}[theorem]{Corollary}

\newtheorem{definition}[theorem]{Definition}

\newtheorem{remark}{Remark}

\newcommand{\C}{\mathbb{C}}
\newcommand{\K}{\mathbb{K}}

\usepackage{xcolor}

\author[Z. Jelonek]{Zbigniew Jelonek}
\author[K. Zaj\k ac]{Karolina  Zaj\k ac}

\address{Z.Jelonek, Instytut Matematyczny, Polska Akademia Nauk, \'Sniadeckich 8, 00-656 Warszawa, Poland.}

\email{najelone@cyf-kr.edu.pl}

\address{K. Zaj\c ac, Instytut Matematyczny UJ, \L ojasiewicza 7, Krak\'ow , Poland.}

\email{karolina.w.zajac@doctoral.uj.edu.pl}

\title[On Algebraic LNE models of complex affine algebraic curves]{On Algebraic LNE models of complex affine algebraic curves}

\begin{document}
	
\begin{abstract}
We give a sufficient condition for an affine complex curve to have an algebraic complex LNE model.
\end{abstract}
			
\maketitle			
			
	\section{Introduction}
			
			There is recently focused interest in metric geometry of algebraic sets, especially in Lipschitz classification problems. 
			Since the result of Birbrair and Mostowski in \cite{birbrairmostowski} has been extended to non-compact case by \cite{OnePoint}, one can classify affine  algebraic sets
			with respect to inner Lipschitz equivalence by means of their LNE models, i.e. sets which are LNE and inner bi-Lipschitz to the initial sets. This means that the inner Lipschitz classification of affine algebraic sets is equivalent to outer Lipschitz classification of a certain class of semialgebraic sets. The tools for outer classification are much stronger and are applied in a restrictive family of sets, thus this approach is very promising.
			
			The essential question is the following: is inner Lipschitz classification of connected algebraic sets equivalent to outer Lipschitz classification of LNE sets which are also algebraic?
			
			In the paper \cite{JMM}  this problem is solved in the case of real curves. The answer is affirmative in the case of real algebraic curves and negative in the case of complex algebraic curves. The aim of this note is to give a sufficient condition for a complex affine curve to have an algebraic LNE model. Our main result is a following:
			
			\vspace{5mm}
			
		\noindent	{\bf Theorem \ref{thm:noComplexmodels}}
			
	\noindent			{\it Let $X\subset \C^n$ be an affine algebraic curve of genus $g$ with $p$ branches at infinity. Assume that apart of simple singularities it has $r$-fold singularities at points $a_1,...,a_k$ of orders $r_1,...,r_k$, respectively. If $$p\ge 2g + 2\sum^k_{i=1}(r_i-1) +1,$$ then $X$ has an algebraic LNE model.}

	\vspace{5mm}

	Let us recall that $X$ has an {\it $r$-fold} singularity of order $r_a$ at the point $a$ if the germ $X_a$ has exactly $r_a$ irreducible components. If all this components are smooth and tangent lines to all branches are linearly independent, we say that this $r$-fold is {\it simple}.

    \vspace{5mm}

As a consequence of the above Theorem, we obtain:

\vspace{5mm}

\noindent {\bf Corollary \ref{genus0}}

\noindent  {\it A smooth affine curve of genus $0$ always has an algebraic LNE model in $\C^3$.}

\vspace{5mm}
We have also the following necessary condition:

\vspace{5mm}

  \noindent	{\bf Theorem \ref{thm:noComplexmodels1}}
{\it Let $X\subset \C^n$ be an affine algebraic curve of genus $g$ with $p$ branches at infinity. Assume that apart of simple singularities it has $r$-fold singularities at points $a_1,\dots,a_k$ of orders $r_1,\dots,r_k$, respectively. If $X$ has an algebraic model, then 
imperatively    $$p^2-3p+2\ge 2g + \sum^k_{i=1}(r_i-1)r_i.$$}
			
			\section{Preliminaries}
			
			\subsection{Basic definition}

			Throughout this subsection let $(M,d)$ and $(N,d')$ be  metric spaces and let $X\subset M$ and $Y\subset N$.

			\begin{definition}
				A map $f:(M,d)\to (N,d')$
				is said to be bi-Lipschitz with constant $L>0$ if for all $x,y\in M$ we have
				$$
				\dfrac{1}{L}d(x,y)\leq d'(f(x),f(y))\leq Ld(x,y).
				$$	
			\end{definition}
			
			A subset $X$ of a metric space can be equipped with the inner (length) metric~$d_X$ given by infimum on lengths of continuous curves connecting given points as well as with the outer metric which is the restriction of the ambient metric to the set (and so we will simply denote the outer metric the same as the ambient metric).

			\begin{definition}
				A  set $X\subset M$ is said to be {\it Lipschitz normally embedded (LNE)} in $M$ if the map $$id:(X,d)\to (X,d_X)$$ is bi-Lipschitz (in other words, the inner metric is equivalent to the outer metric).
			\end{definition}
			
			This definition in particular implies that the inner metric on an LNE set $X$ is a~well-defined distance function and any two points in $X$ can be joined by a curve of finite length.

			\begin{definition}
				We say that $f:X\to Y$ is \emph{outer bi-Lipschitz} if it is bi-Lipschitz with respect to the outer metrics on $X$ and $Y$ and we say $f$ is \emph{inner bi-Lipschitz} if it is bi-Lipschitz with respect to the inner metrics.
			\end{definition}
			
			Note that  if two sets $X$ and $Y$ are outer bi-Lipschitz, then they are also inner bi-Lipschitz.
			In particular an outer bi-Lipschitz image of an LNE set is LNE.

			\begin{definition}[LNE models]
				Let  $f:X\to N$ be a continuous map. We say that $Y$ is an LNE model of $X$ in $N$ when the set $Y:=f(X)$ is LNE and the map $f:(X,d_X)\to (Y,d_Y)$ is inner bi-Lipschitz.
			\end{definition}

\subsection{Important results}

We will use the crucial method of extending polynomial isomorphisms.

\begin{theorem}[\cite{jel},\cite{jel1}]\label{J}
Let $X$ be a nonsingular (not necessarily connected) 
closed subvariety of $\mathbb{C}^n$
of dimension (not necessarily pure) $k$. Let us consider a regular embedding $f = (p_1, ..., p_n) : X \to \mathbb{C}^n$. If $n \ge 2k + 2$, then there exists a tame polynomial isomorphism $F =
(P_1, ..., P_n) : \mathbb{C}^n \to \mathbb{C}^n$ such that $F_{|X} = f.$
\end{theorem}

\begin{theorem} [\cite{JMM2}, \cite{buz}]\label{wektory}
Let $\K$ be an infinite field. Let  $n\ge 2$ and $a_1,\ldots, a_k; $ $b_1,\ldots , b_k$ be two sequences of different points from $\K^n.$ Let $L_1,\ldots, L_k\in SL(n,\K)$ be linear isomorphisms with 
jacobian one. 
Then there is a polynomial automorphism $\Phi: \K^n\to \K^n$ with jacobian one such that $\Phi(a_i)=b_i$ and $d_{a_i} \Phi= L_i$ for $i=1,\ldots, k.$
\end{theorem}

The following theorem will also be useful:

			\begin{theorem}[\cite{cgm}]\label{cgm}  Let $X$ be a connected projective curve of $\C\mathbb{P}^n$ of degree~$d$.  Assume no component of~$X$ is contained in the hyperplane at infinity $H_\infty$. 
				The affine algebraic curve $X_a:=X\setminus H_\infty$ is LNE in $\C^n$ if and only if the following conditions are
				satisfied: 
				
				(0) $X_a$ is connected; 
				
				(i) $X_a$ is locally LNE at each of its singular points; 
				
				(ii) The curve X intersects $H_\infty$ transversally, i.e., in exactly $d=\deg X$ different points.
			\end{theorem}

    \begin{remark}    {\rm It is well known (\cite{np}), that a curve $X$ is locally LNE at a point $a$ if it has a normal $r$ fold singularity (with $r\ge 1$) at $a$, i.e., it is an union of $r$ smooth branches, which intersect pairwise transversally.}
\end{remark}

    Here we prove theorem about topological invariance of the inner metric of a complex curve. Our proof is in fact a modification
    of a local proof of Pichon-Neumann. It is worth to note that a version for $n=2$ was proved in \cite{fs}.
	
	\begin{theorem} \label{top} 
    Let $X_1, X_2 \subset\C^n$ be two complex irreducible algebraic curves. Then, the
following statements are equivalent:

(1) $X_1$ and $X_2$ are homeomorphic;

(2) $X_1$ and $X_2$ are inner lipeomorphic, i.e. bi-Lipschitz homeomorphic with respect
to the inner metrics.
\end{theorem}

\begin{proof}
The implication $(\Rightarrow)$ is trivial. For $(\Leftarrow)$. We decompose any connected affine algebraic curve $X \subset \mathbb{C}^n$ into standard metric models determined solely by topological invariants.

\medskip
\noindent\textbf{Step 1: Geometric decomposition.} 
$X$ has finitely many singular points $\operatorname{Sing}(X) = \{p_1, \dots, p_s\}$ and finitely many ends at infinity $E_1, \dots, E_k$ (corresponding to branches of $\overline{X} \subset \mathbb{CP}^n$ along $H_\infty$).
\begin{itemize}
    \item \textbf{Singularities:} Surround each $p_i$ with a small disjoint ball $B(p_i, \varepsilon_i)$ whose boundary $S(p_i, \varepsilon_i)$ intersects $X$ transversally. Each $U_i = X \cap \overline{B(p_i, \varepsilon_i)}$ splits into $r_i \ge 1$ irreducible local branches $U_{i,1}, \dots, U_{i,r_i}$ meeting only at $p_i$.
    \item \textbf{Ends:} For sufficiently large $R > 0$ (enclosing all $\overline{B(p_i, \varepsilon_i)}$ with $S(0,R)$ transversal to $X$), $X \setminus B(0,R)$ decomposes into $k$ disjoint ends $E_1 \sqcup \dots \sqcup E_k$.
    \item \textbf{Compact core:} Removing the interiors of all $U_i$ and $E_j$ leaves a compact, smooth, orientable surface $K$ with $\sum_{i=1}^s r_i + k$ smooth boundary circles.
\end{itemize}

\medskip
\noindent\textbf{Step 2: Straightening local branches.} 
Translating $p_i$ to $0$, Puiseux's theorem yields a parametrization $\gamma \colon \overline{\mathbb{D}}_\delta \to U_{i,l}$ of the form $\gamma(t) = (t^m, \sum_{\nu \ge m} a_{2,\nu} t^\nu, \dots)$, where $m \ge 1$ is the multiplicity. For small $\delta > 0$, $\|\gamma'(t)\| \asymp |t|^{m-1}$ on $\overline{\mathbb{D}}_\delta \setminus \{0\}$, where $a \asymp b$ means that there exists positive constant $c,C$ such that $cb<a<Cb$. Composing $\gamma^{-1}$ with the polar map $\Phi(r e^{i\theta}) = (r/\delta)^m e^{i\theta}$ gives a diffeomorphism $h_{i,l} = \Phi \circ \gamma^{-1} \colon U_{i,l} \setminus \{p_i\} \to \overline{\mathbb{D}} \setminus \{0\}$ satisfying
\[
ds_{\overline{\mathbb{D}}}^2 = \frac{r^{2(m-1)}}{\delta^{2m}}\bigl(m^2 dr^2 + r^2 d\theta^2\bigr) \asymp |t|^{2(m-1)}|dt|^2 \asymp ds_{U_{i,l}}^2.
\]
Since a point has $1$-dimensional measure zero, $h_{i,l}$ extends across $p_i$ to an inner bi-Lipschitz homeomorphism $U_{i,l} \to \overline{\mathbb{D}}$.

\medskip
\noindent\textbf{Step 3: Straightening ends at infinity.} 
Applying Puiseux's theorem to the branch at $H_\infty$ corresponding to $E_j$ with $z = 1/t$, we obtain a parametrization $\psi \colon \{t \in \mathbb{C} : |t| \ge R_0\} \to E_j$ given by Laurent series $\psi(t) = (t^m, \sum_{\nu \le m} b_{2,\nu} t^\nu, \dots)$ with pole order $m \ge 1$. For large $R_0$, $\|\psi'(t)\| \asymp |t|^{m-1}$. Using the polar map $\Psi(r e^{i\theta}) = (r/R_0)^m e^{i\theta}$, the exact same length element estimate shows that $g_j = \Psi \circ \psi^{-1} \colon E_j \to \mathbb{C} \setminus \mathbb{D}$ is an inner bi-Lipschitz homeomorphism.

\medskip
\noindent\textbf{Step 4: Gluing on the compact core.} 
A homeomorphism $f \colon X_1 \to X_2$ preserves the genus $g$, the number of ends $k$, and (by invariance of domain) induces a bijection between multibranch singularities ($r_i \ge 2$) preserving branch counts $r_i$. (Unibranch cusps with $r_i = 1$ are topological disks and are treated uniformly via Step 2.) Thus, the compact cores $K_1 \subset X_1$ and $K_2 \subset X_2$ have the same genus $g$ and $N = \sum_{i=1}^s r_i + k$ boundary components.

By the classification of compact surfaces with boundary, there exists a diffeomorphism $F_K \colon K_1 \to K_2$ matching $(h_{i,l}^{(2)})^{-1} \circ h_{i,l}^{(1)}$ and $(g_j^{(2)})^{-1} \circ g_j^{(1)}$ on $\partial K_1$ (as orientation-preserving diffeomorphisms of $S^1$ are isotopic). By compactness of $K_1$, $F_K$ is inner bi-Lipschitz. Gluing $F_K$ with the maps on $U_{i,l}^{(1)}$ and $E_j^{(1)}$ defines a global homeomorphism $F \colon X_1 \to X_2$. Subdividing any rectifiable path into segments within individual pieces and taking the maximum Lipschitz constant $L$ across all finitely many pieces yields
\[
L^{-1} d_{X_1}(x,y) \le d_{X_2}\bigl(F(x), F(y)\bigr) \le L d_{X_1}(x,y) \quad \forall x,y \in X_1. \qedhere
\]
\end{proof}

	\section{Complex  algebraic LNE models}\label{sec:algebraicModels}

	Now we discuss LNE models of complex algebraic curves in $\C^n.$ Clearly, by \cite{birbrairmostowski} and \cite{OnePoint} a complex algebraic curve has a semialgebraic LNE model. Yet, \cite{JMM}  shows that  complex algebraic curves may not admit complex algebraic LNE models. This is unlike in the real setting discussed in \cite{JMM}. Here we give a sufficient condition for a~complex affine curve to have an algebraic LNE model.
            
    We start with the following:

    \begin{theorem}[\cite{mum}]
        Let $X$ be a projective curve with an arithmetic genus $p_a(X)$. Let $D$ be a Cartier divisor on $X$.
        If $\deg D\ge 2p_a(X)+1,$ then $D$ is a very ample divisor.
    \end{theorem}

    \begin{corollary}\label{ample}
        Let $X$ be a projective curve of genus $g$ with only simple $r$-fold singularities at points $a_1,..,a_k$  of orders $r_1,\dots,r_k$, respectively.
        If $D$ is a Cartier divisor on $X$ such that $\deg D\ge 2g+2\sum^k_{i=1} (r_1-1)+1,$ then $D$ is very ample.
    \end{corollary}

    \begin{proof}
        We use here the Serre formula $p_a(X)=g+\sum_{a\in Sing(X)} \delta_a.$ It is not difficult to check that in our case $\delta_{a_i}=r_i-1.$
        \end{proof}

	\begin{theorem}\label{thm:noComplexmodels}
		Let $X\subset \C^n$ be an affine algebraic curve of genus $g$ with $p$ branches at infinity. Assume that apart of simple singularities it has $r$-fold singularities at points $a_1,\dots,a_k$ of orders $r_1,\dots,r_k$, respectively. If $$p\ge 2g + 2\sum^k_{i=1}(r_i-1) +1,$$ then $X$ has an algebraic LNE model in $\C^3$.
\end{theorem}

	\begin{proof}
		By Theorem \ref{cgm} and Theorem \ref{top} it is enough to construct an affine curve $Y\subset \C^N$ of genus $g$ for some $n\in \mathbb{N}$, which has exactly $k$ simple $r$-folds of orders $r_1,\dots,r_k$, respectively, and which is transversal to infinity in exactly $p$ points. We do this by induction with respect to the number $s=\sum^k_{i=1} r_i$. 

        For $s=0$, let $X'$ be a normalization of a projective closure of $X$. Let points $b_1,\dots,b_p\in X'$ correspond to branches at infinity of the curve $X$. Consider the divisor $D=\sum b_i.$ By Corollary \ref{ample} and by our assumption, the divisor $D$ is very ample. Hence there exists an embedding of $X'$ in some projective space such that all images of points $b_1,\dots, b_p$ belong to one hyperplane. Taking this hyperplane as a hyperplane at infinity, we obtain a smooth model of the curve $X$.

        Now assume that we have constructed a model $X_s$ with only simple $r$-folds of orders $r'_i\le r_i$ for which $s=\sum^k_{i=1} r'_i < \sum^k_{i=1} r_i.$ We show how to construct such a~model for $s'=s+1.$

        Let $M=max \{r_1,\dots, r_k, 3\}.$ Put $l<k$. Let us assume that $X_s\subset \C^N$ for some $N>2(M+p+1)+2$ and the model $X_s$ has simple $r$-folds of orders $r_1,..,r_l$ at points $a_1,\dots,a_l$, respectively, but $r_{l+1}$ fold is either not yet constructed at all or $r'_{l+1}<r_{l+1}.$ If $r'_{l+1}$ fold in not yet constructed, we choose $a_{l+1}$ to be any point different from $a_1,\dots,a_l.$ 

        Using Theorem \ref{J} , we can assume that all points $a_1,\dots,a_l,a_{l+1}$ belong to $\C^M\times\{0\},$ but the whole curve $X$ is not contained in $ \C^M\times\{0\}$. Moreover, we can assume that in $\C^M\times\C^{p+1}\times \{0\}$ there exist pairwise distinct points $b_1,\dots,b_{p+1} \in X$, which are not in $\C^m\times \{0\}$ but which determine all coordinate axes in $\{0\}\times \C^{p+1}\times \{0\}.$ Again using, Theorem \ref{J} , we can assume that the all points of $X$ in $\C^M\times\C^{p+1}\times \{0\}$ are exactly points in $\C^m\times \{0\}$ together with these extra points $b_1,\dots,b_{p+1}$.

        Using Theorem \ref{wektory}, we can assume that tangent vectors to the model $X_s$ in points $b_1,\dots,b_{p+1}$ are parallel to $\C^M\times \{0\}$ and linearly independent with tangent vectors of branches of $X_s$ at the point $a_{l+1}.$ We can also assume that in $\C^M\times \{0\}$ there exists a simple point $c$ different from all $a_1,\dots,a_k$.

        Now consider $p+1$ projections $\pi_0, \pi_1, \dots, \pi_{p+1}$ along lines $a_{l+1}, b_1,\dots,b_{p+1}$ into hyperplanes of the form $\C\times \C \dots \times\{0\}\times \C \dots \times \C$ (with \{0\} on $(M+i+1)$-th place for $\pi_i$). Since the curve $X_s$ has at most $p$ points at infinity, at least one of these projections is proper. Moreover, it transforms the simple point $c$ into the same point $c$ with multiplicity one. It means that this mapping is birational. It preserves all points $a_i$ unchanged. Hence it can glue only finite points different than $a_i$. Using a series of blowing up's if necessary, we can modify the curve in such a~way that it has only simple $r$-folds in $a_1,\dots,a_{l+1}$. Using again Corollary \ref{ample}, we obtain a new model $X_s'$ with $s'>s.$

        We finish the proof by induction. Finally, it is enough to use a generic projection onto $\C^3$, which is bi-lipschitz by \cite{bfj}.  
    \end{proof}

	\begin{corollary}\label{genus0}
        A smooth affine curve of genus $0$ always has an algebraic LNE model in $\C^3.$
	\end{corollary}

 Now we are looking for necessary conditions:

    \begin{theorem}\label{thm:noComplexmodels1}
		Let $X\subset \C^n$ be an affine algebraic curve of genus $g$ with $p$ branches at infinity. Assume that apart of simple singularities it has $r$-fold singularities at points $a_1,\dots,a_k$ of orders $r_1,\dots,r_k$, respectively. If $X$ has an algebraic model, then 
imperatively    $$p^2-3p+2\ge 2g + \sum^k_{i=1}(r_i-1)r_i.$$
	\end{theorem}

\begin{proof}
    Assume that such an algebraic model $Y\subset\C^m$ exists. Let $\pi: Y\to\C^2$ be a generic projection. It has exactly $p$ smooth ponts at infinity at which is transversal to the line at infinity. Hence it has degree $p.$ This measn that $p_a(Y)=(p-1)(p-2)/2.$ Moreover it has $k$ normal r-folds of arder $r_1,..., r_k$ and maybe $s$ folds produced by a projection. Hence by Serre formula
    $$(p-1)(p-2)/2=g+\sum^k_{i=1}(r_i-1)r_i/2+s\ge g+\sum^k_{i=1}(r_i-1)r_i/2.$$
\end{proof}

\section{Acknowledgments}
 The first named author was partially supported by the grant of Narodowe Centrum Nauki, number 2024/55/B/ST1/01412.

\end{document}